\documentclass{esannV2}
\usepackage[english]{babel}
\usepackage{graphicx}
\usepackage[utf8]{inputenc}
\usepackage[fleqn]{amsmath}
\usepackage{amssymb,array}
\usepackage{subfigure}
\usepackage{array}
\usepackage{multirow}
\usepackage{url}
\usepackage{enumitem}

\begin{document}
\title{Activity Date Estimation in Timestamped Interaction Networks}

\author{Fabrice Rossi and Pierre Latouche
\vspace{.3cm}\\
SAMM EA 4543, Université Paris 1 Panthéon-Sorbonne\\
90, rue de Tolbiac, 75634 Paris cedex 13, France\\
}

\maketitle

\begin{abstract}
We propose in this paper a new generative model for graphs that uses a latent
space approach to explain timestamped interactions. The model is designed to
provide global estimates of activity dates in historical networks where only
the interaction dates between agents are known with reasonable
precision. Experimental results show that the model provides better results
than local averages in dense enough networks.
\end{abstract}

\section{Introduction}
In this paper, we study interactions between agents that are recorded on a time
scale larger than the expected lifespan of the agents. A typical instance of
such interactions are property ownership recordings in which a house or a land
exists for a very long time period and passes from owner to owner, outliving
them. Ownerships are generally recorded with a lot of details in well kept
archives, while the lives of the owners are generally known with much less
details. 

In
\cite{bouletjouveetal2008neurocomputing,rossivilla-vialaneixetal2011exploration-large}
for instance, the primary information source consists in notarial acts
recording different forms of ownerships of lands and related objects (according to
the French feudal laws) during a long time period (up to 300 years for the digitized
version). On the one hand, most of the notarial acts have a proper date,
precise at least at the year level, while, on the other hand, little is known
about the tenants (a.k.a. the ``owners'') involved in the acts, apart from
their names. It seems therefore interesting to infer information about the
tenants from the acts, namely to estimate a living period for each tenant
based on the acts in which he/she is involved.

More generally, we consider a graph whose vertices represent agents and
whose edges represent interactions between those agents. We consider here
simple graphs, but the approach generalizes immediately to multi-graphs in
which several edges can link two agents. Each interaction is timestamped and
our goal is to estimate a central time stamp for each agent in such a way that
interaction dates are compatible with the time stamps of the agents and with
expert knowledge on the expected life span of the agents. The main difficulty
of this task comes from inconsistencies observed in real world historical
data: due to name ambiguities, associations between agents and interactions
are sometimes incorrect. In network parlance this corresponds to some rewiring
of the graph: while we should get a connection between $a$ and $b$, a naming
ambiguity between vertex $b$ and $c$ assigns wrongly this interaction to $a$
and $c$. 

We propose a solution based on a generative model inspired by the latent space
model of \cite{hoffelal2002latent-space}: given the interaction dates, the
model generates interaction networks that fulfill the compatibility constraints
exposed above. Note that the proposed setting is quite different from
classical temporal graph modeling (see
e.g. \cite{goldenbergzhengetal2009survey,xingfuetal2010StateSpace}) where the primary goal generally
consists in understanding the evolution of the structure of the network
through time. 

The rest of the paper first introduces the generative model as well as the maximum
likelihood estimation strategy. It then summarizes experimental results on
simulated data.

\section{A Generative Model}
We observe an undirected graph $G$ characterized by a vertex set $V=\{1,\ldots,n\}$ and
a binary adjacency matrix $A$. When $A_{ij}=1$, that is when node $i$ is
connected to node $j$, we are given an associated interaction date
specified as a positive real number $D_{ij}$. 

We consider a generative model for $(A,D)$ based on latent activity date
variables. More precisely, each vertex $i$ is associated to a positive
(unobserved) real number $Z_i$ which summarizes the activity period of said
vertex. Then, we assume that the probability of having a connection between $i$
and $j$ is linked to the temporal distance $|Z_i-Z_j|$. We assume also that
knowing $Z=(Z_i)_{1\leq i\leq n}$, the $A_{ij}$ are independent. Finally, when $i$ and $j$
are connected, we assume that their interaction date is randomly distributed
between $Z_i$ and $Z_j$ (independently of all other variables). In more
technical terms, the conditional independence assumptions lead to the
following generative model, where $\theta$ denotes numerical parameters:
\begin{multline}
  \label{eq:generative}
p(A,D|Z,\theta)=\prod_{i\neq j,A_{ij}=0}P(A_{ij}=0|Z_i,Z_j,\theta)\\
\times\prod_{i\neq j,A_{ij}=1}p(D_{ij}|A_{ij}=1,Z_i,Z_j,\theta)P(A_{ij}=1|Z_i,Z_j,\theta).
\end{multline}

\subsection{A specific model}
We specialize now the generic form of equation
\eqref{eq:generative}. Inspired by \cite{hoffelal2002latent-space}, we use a
logistic regression model for the connection probabilities, that is
\begin{equation}
  \label{eq:connection}
\log\frac{P(A_{ij}=1|Z_i,Z_j,\alpha,\beta)}{P(A_{ij}=0|Z_i,Z_j,\alpha,\beta)}=\alpha-\beta(Z_i-Z_j)^2,
\end{equation}
while the interaction date $D_{ij}$ is simply modelled with a Gaussian
distribution around $\frac{Z_i+Z_j}{2}$, that is
\begin{equation}
  \label{eq:dates}
D_{ij}|Z_i,Z_j,\sigma\sim \mathcal{N}\left(\frac{Z_i+Z_j}{2},\sigma^2\right).
\end{equation}
Then, up to constants, the log-likelihood of the data is given by
\begin{multline}
  \label{eq:log_likelihood}
L(A,D|Z,\sigma,\alpha,\beta)=\sum_{i\neq j,A_{ij}=1}\left(-\log\sigma-\frac{1}{2\sigma^2}\left(D_{ij}-\frac{Z_i+Z_j}{2}\right)^2\right)\\
+\sum_{i\neq j}\left(A_{i,j}(\alpha-\beta(Z_i-Z_j)^2)-\log\left(1+e^{\alpha-\beta(Z_i-Z_j)^2}\right)\right).
\end{multline}
Connection probabilities are not identical to the ones used in
\cite{hoffelal2002latent-space} for two reasons. Firstly, we use a quadratic
term $(Z_i-Z_j)^2$ rather than the original absolute value $|Z_i-Z_j|$ to
avoid numerical instabilities linked to the non differentiability of the
latter. Secondly, we add a $\beta$ parameter to compensate for the relatively
large values found in real world historical networks for $(Z_i-Z_j)^2$ which
can be of the order 2500. In \cite{hoffelal2002latent-space}, the absence of
the first term in equation \eqref{eq:log_likelihood} allows for free scaling
effects of the $Z_i$, something that is not possible here.

\subsection{Estimation}
We use a maximum likelihood approach implemented via a gradient descent based algorithm. A
natural initialization for $Z$ is provided by
$\hat{Z}_i=\frac{\sum_{j,A_{ij}=1}D_{ij}}{\sum_{j}A_{ij}}$, where $\hat{Z}_i$ takes the
average value of the dates of the outgoing/incoming edges. This corresponds
to an estimation of the activity dates based only on local information. The initial
values of the other parameters are chosen as follows. As $\sigma$ models the
life span of actors, we use $50$ as a starting point, allowing interactions to
happen in a very large two hundred years interval centered in the average activity
date of two actors. The $\alpha$ parameter is initialized such that the
connection probability equals the observed network density when all the $Z_i$
are equals. Finally, $\beta$ is set to a value that reduces the connection
probability to almost zero ($10^{-6}$) when the temporal distance between two
actors is above one hundred years. 

\section{Experimental evaluation}
As we do not have access yet to large historical data sets with reliable activity
dates, we focus on simulated data to evaluate the model we propose and to understand its
strengths and limitations.

\subsection{Ideal situation}
First, we study the ideal situation in which networks are simulated according
to our model. The goal is then to recover the activity dates used to generate
the data. We measure the quality of the recovery using the mean square error
(MSE) between the real $Z$ and the estimated one. A network is generated as
follows:
\begin{enumerate}[noitemsep]
\item $n=100$ and the $Z_i$ are uniformly distributed in $[1200,1400]$;
\item for a given maximal target density $d$ in $[0.1,0.5]$, $\alpha$ is set
  to $-\log\left(\frac{1}{d}-1\right)$. This sets the probability of
  $A_{ij}=1$ to $d$ when $Z_i=Z_j$; 
\item the expected life span is set to 80. Accordingly, $\beta$ is set to 
$\frac{\log\left(\frac{1}{\varepsilon}-1\right)+\alpha}{80^2}$,
  where $\varepsilon$ is the target probability to have $A_{ij}=1$ when
  $|Z_i-Z_j|=80$. We use $\varepsilon=10^{-6}$;
\item $\sigma$ is set to $20$ (that is to one fourth of the life span);
\item given $\alpha$, $\beta$, $\sigma$ and $Z$, $A$ and $D$ are generated
  according to the model;
\item finally, we keep only the largest connected component of the obtained
  graph. We discard graphs in which there are less edges than the number of
  parameters in the model (that is three added to the number of vertices). 
\end{enumerate}

\begin{figure}[htbp]
\mbox{}\vskip -1.5em
  \centering
\includegraphics[width=0.85\linewidth]{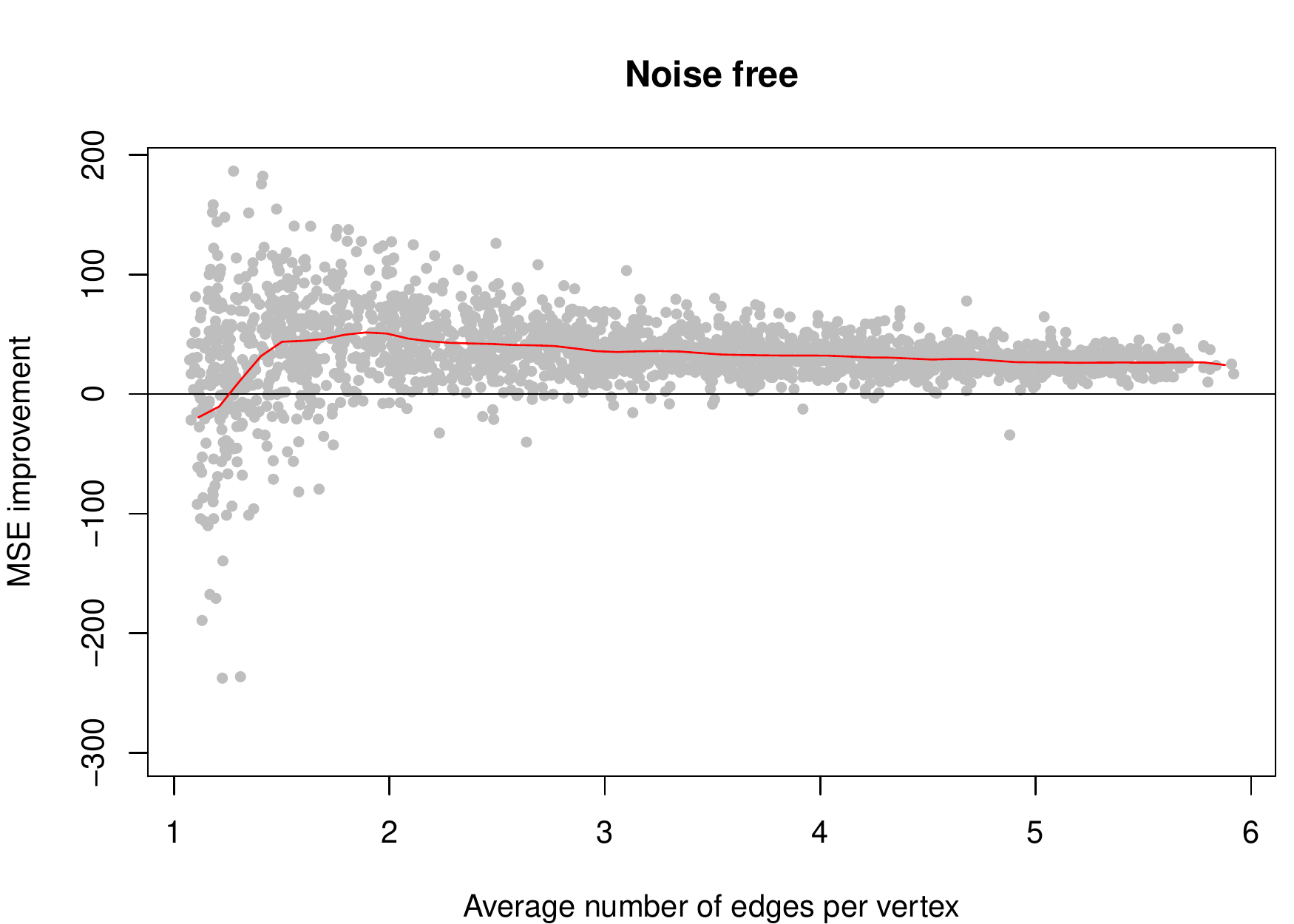}  
  \caption{\small Ideal situation: each dot gives the improvement in MSE over
    local averages when using the proposed model estimates as a function of
    the number of edges per node. On a total of 2152 networks, 3 where
    excluded from the figure because of very large negative values of the improvement (down to -2200) due to convergence issues. Those networks had below 1.27 edges per vertex.}
  \label{fig:ideal}
\mbox{}\vskip -1.5em
\end{figure}

Results are summarized by Figure \ref{fig:ideal} which gives the improvement
in MSE obtained by using the $Z$ estimates of the model compared to the local
averages $\hat{Z}_i$ defined above. The
superimposed curve is a kernel based estimate of the relation 
between the average number of edges per vertex and the improvement in
MSE. According to this estimate, the average improvement reaches a positive
value above 1.31 edges per vertex. When this number is above 2, the
improvement over local estimates is almost always positive. 

\subsection{Misspecification}
\begin{figure}[htbp]
\mbox{}\vskip -1.5em
  \centering
\includegraphics[width=0.85\linewidth]{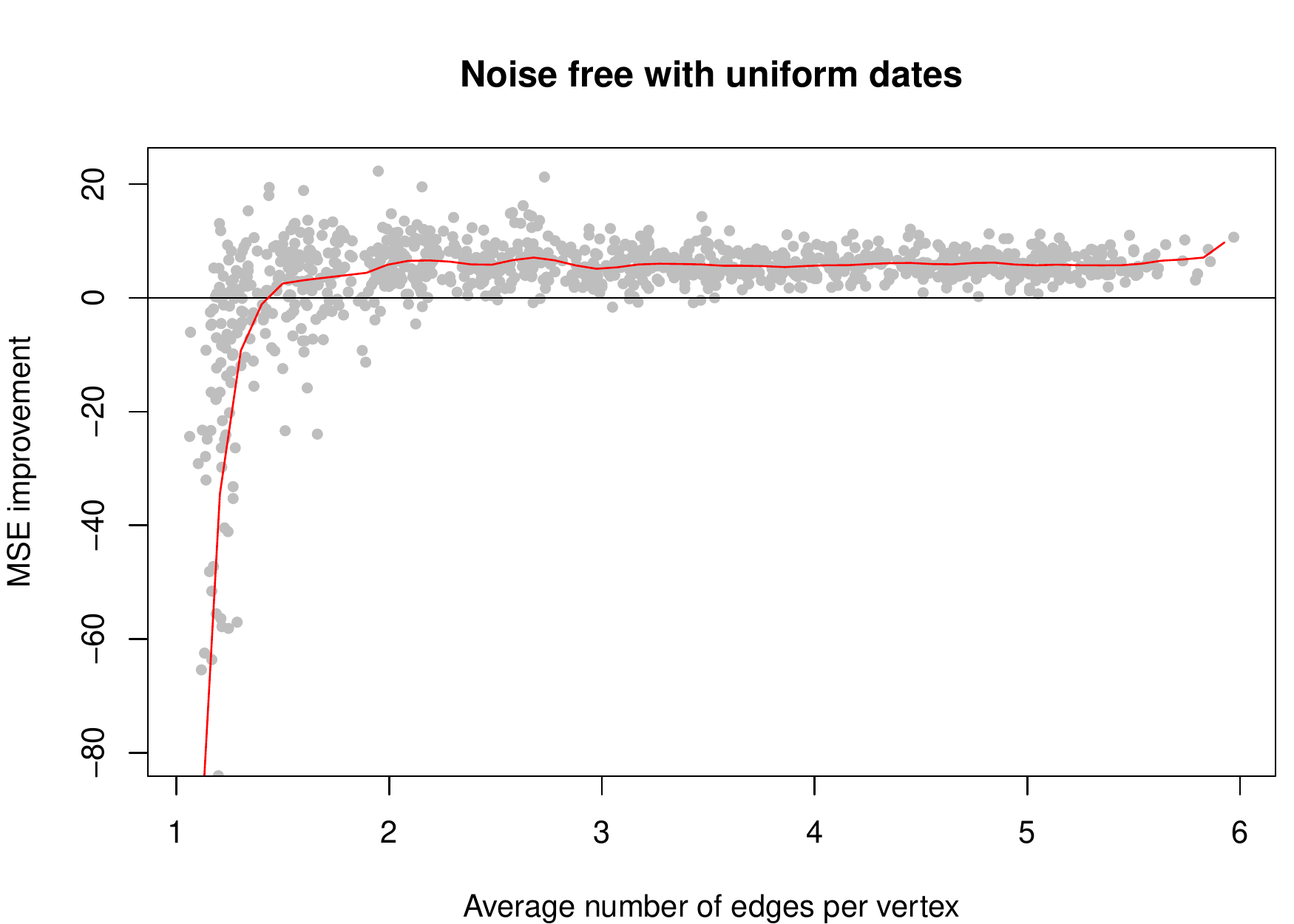}  
\caption{\small Misspecification with uniform connection dates. 10 graphs out
  of 1079 with an average number of edges per vertex below 1.21 and bad
  estimates were removed from the figure to keep it readable.}
  \label{fig:uniform}
\mbox{}\vskip -1.5em
\end{figure}

We also tested the model under mis-specification by replacing the
Gaussian distribution for dates by a uniform distribution between $Z_i$ and
$Z_j$ for connected vertices. This introduces some form of
heteroscedasticity. Results displayed on Figure \ref{fig:uniform} show a good
resistance of the model to misspecification when the number of edges per
vertex if above 2. 

\subsection{Rewiring}
\begin{figure}[htbp]
\mbox{}\vskip -1em
  \centering
\includegraphics[width=0.85\linewidth]{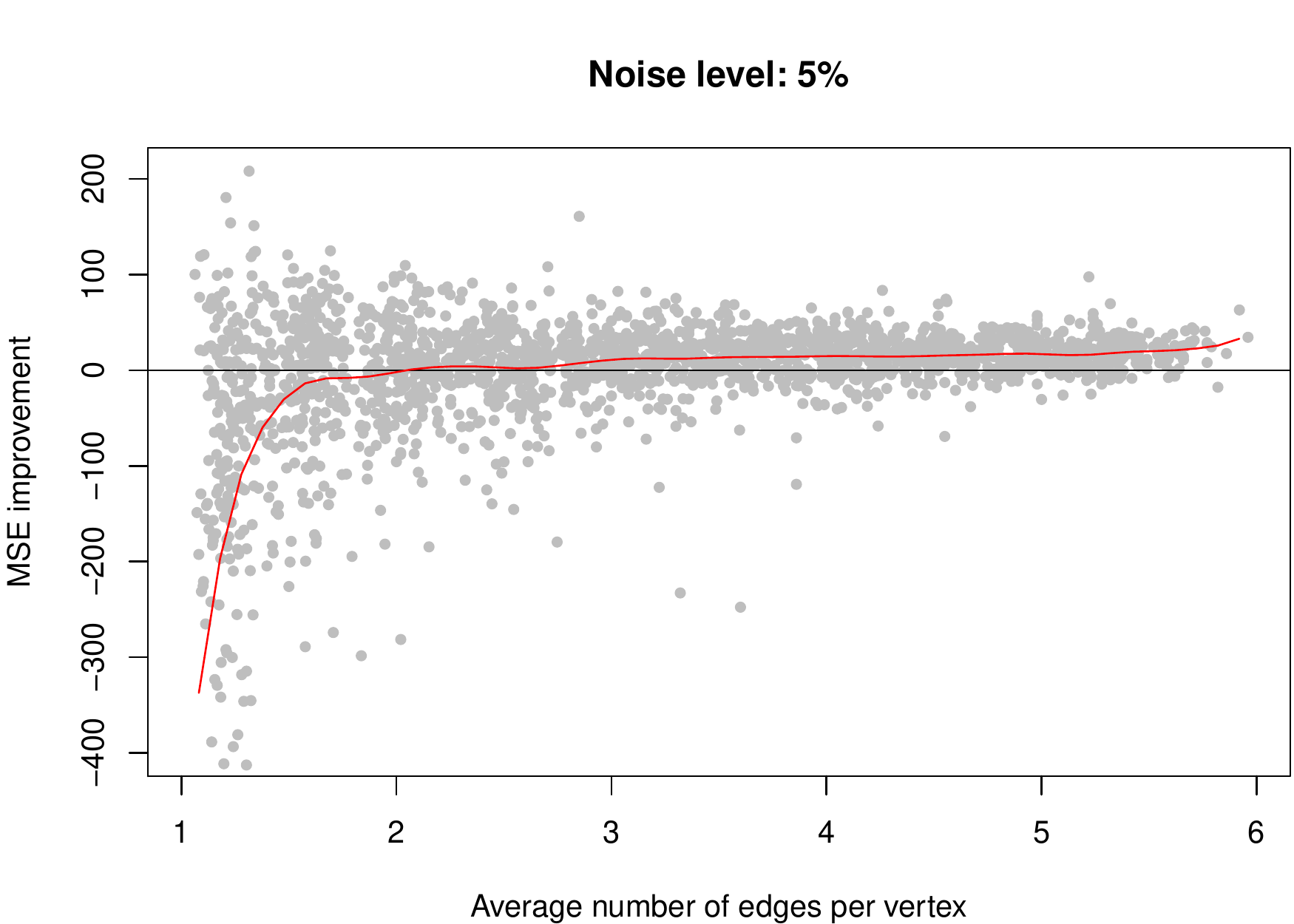}  
\caption{\small Improvements under rewiring noise. 54 graphs out
  of 2159 with an average number of edges per vertex below 1.79 and very bad
  estimates are removed from the figure to keep it readable.}
  \label{fig:noise}
\mbox{}\vskip -1em
\end{figure}

Finally, we study the robustness of the model with respect to the rewiring
issue exposed in the introduction. Networks are first generated according to
the model and then a certain number of edges are randomly rewired by moving
one of the end points to a randomly selected vertex while keeping the original
date. In the case of a very low noise (1\% of rewired edges), almost no effect
on the improvements are observed (results not shown here). Figure
\ref{fig:noise} shows results for a more important noise (5\% of rewired
edges). As expected, the model, while showing robustness, is impaired by the
``false'' information attached to rewired edges. According to the kernel
estimator, at least 2.1 edges per vertex are needed to reach equal
performances between the local average and the model estimates, while above 3,
the model outperforms the local estimates.  

\section{Conclusion}
Results on simulated data are very satisfactory: above an average number of two
edges per vertex, the estimates provided by the model are closer to the ground
truth than local averages, even under two forms of misspecification (uniform
date distribution and edge rewiring). While the estimator exhibits large
variability and can give quite bad results, this happens only under a low
number of edges per vertex. While a direct numerical evaluation of the method
on real world historical data is impossible because of the lack a reliable
activity dates on large interaction databases, we are working with historians
on qualitative assessment of the results based on dates inferred from well
known figures such as prominent land lords. 


\begin{footnotesize}
\bibliographystyle{abbrv}
\bibliography{bibliography}
\end{footnotesize}

\end{document}